\documentclass[a4paper,11pt,reqno]{amsart}
 \usepackage{amssymb,amsmath,amstext,amsthm,amsfonts}
\usepackage{euscript}

\newtheorem{theorem}{Theorem}

\newtheorem{lemma}{Lemma}
\newtheorem{proposition}{Proposition}

\theoremstyle{definition}

\newtheorem{remark}{Remark}

\newcommand{\ep}{\varepsilon}

\newcommand{\bR}{{\bf R}}

\newcommand{\bZ}{{\bf Z}}
\newcommand{\bN}{{\bf N}}

\newcommand{\bS}{{\bf S}}
\newcommand{\bT}{{\bf T}}

\newcommand{\cA}{{\mathcal A}}
\newcommand{\cB}{{\mathcal B}}
\newcommand{\cC}{{\mathcal C}}

\newcommand{\cW}{{\mathcal W}}
\newcommand{\cM}{{\mathcal M}}
\newcommand{\cU}{{\mathcal U}}

\newcommand{\cP}{{\mathcal P}}
\newcommand{\cE}{{\mathcal E}}

\newcommand{\de}{{\delta}}
\newcommand{\be}{{\beta}}
\newcommand{\al}{{\alpha}}

\renewcommand{\epsilon}{\varepsilon}
\newcommand{\te}{\theta}

\DeclareMathOperator{\inte}{int}
\DeclareMathOperator{\diff}{Diff}
\DeclareMathOperator{\diam}{diam}

\DeclareMathOperator{\dist}{dist}
\DeclareMathOperator{\supp}{supp}

\DeclareMathOperator{\htop}{h_{\rm top}}

\title[Semicontinuity of entropy] {Semicontinuity of entropy, existence of
  equilibrium states and continuity of physical measures}

\begin{thanks} {The author was partially supported by CNPq
    (Brazil), grant BPD/16082/2004 (FCT-Portugal) and
    POCI/MAT/61237/2004 (FCT/POCI 2010
    Portugal/FEDER-EU). Part of this work was done while
    enjoying a post-doc period at PUC (Rio de Janeiro) and
    IMPA whose facilities and research atmosphere are
    gratefully acknowledged.}
\end{thanks}

\author[V\'\i tor Ara\'ujo]{}
\email{vitor.araujo@im.ufrj.br}
 \email{vdaraujo@fc.up.pt}



\begin{document}

\maketitle

\centerline{\scshape V\'\i tor Ara\'ujo }
\medskip

{\footnotesize
\centerline{Instituto de Matem\'atica, Universidade Federal
  do Rio de Janeiro}
\centerline{C. P. 68.530,
  21.945-970 Rio de Janeiro, RJ-Brazil}
\centerline{and}
\centerline{ Centro de Matem\'atica da Universidade do Porto }
\centerline{ Rua do Campo Alegre 687, 4169-007 Porto, Portugal  }
}


\begin{abstract}
    We obtain  results on
    existence and continuity of physical measures through
    equilibrium states and apply these to non-uniformly
    expanding transformations on compact manifolds with
    non-flat critical sets, deducing sufficient conditions
    for continuity of physical measures and, for local
    diffeomorphisms, necessary and sufficient conditions for
    stochastic stability.
    In particular we show that, under certain conditions,
    stochastically robust non-uniform expansion implies
    existence and continuous variation of physical measures.

\end{abstract}

\section{Introduction and statement of results}
\label{sec:introduction}

The statistical viewpoint on Dynamical Systems is one of the
cornerstones of most recent developments in dynamics. Given
a map $f_0$ from a manifold $M$ into itself, a central
concept is that of \emph{physical measure}, a
$f_0$-invariant probability measure $\mu$ whose
\emph{ergodic basin}
\[
B(\mu)=\big\{x\in M:
\frac1n\sum_{j=0}^{n-1}\varphi(f_0^j(x))\to\int\varphi\,
d\mu\mbox{  for all continuous  } \varphi: M\to\bR\big\}
\]
has positive \emph{volume} or \emph{Lebesgue measure}, which
we write $m$ and take as the measure associated with any
non-vanishing volume form on $M$.

This kind of measures provides asymptotic information on a
set of trajectories that one hopes is large enough to be
observable in real-world models.

Here we present recent developments on the relation between
the existence of physical measures and of equilibrium
  states for smooth endomorphisms $f_0:M\to M$ of a compact
boundaryless finite dimensional Riemannian manifold.
We obtain some results of existence and continuity of
physical measures through equilibrium states and apply these
to non-uniformly expanding transformations on compact
manifolds, deducing sufficient conditions for continuity
of physical measures and necessary and sufficient conditions
for stochastic stability.

The stability of physical measures under small variations of
the map
allows for small errors along orbits not to disturb too much
the long term behavior, as measured by the most basic
statistical data provided by asymptotic time averages of
continuous functions along orbits.  In principle when
considering practical systems we cannot avoid external
noise, so every realistic mathematical model should exhibit
these stability features to be able to cope with unavoidable
uncertainty about the ``correct'' parameter values, observed
initial states and even the specific mathematical
formulation involved.


\subsection{Pressure and Equilibrium States}
\label{sec:press-equil-stat}

The physical measures  are related to
equilibrium states of a certain potential function. Let
$\phi : M \rightarrow \bR$ be a continuous function.
Then a $f_0$-invariant probability measure $\mu$ is a
\emph{equilibrium state for the potential $\phi$} if
\[
P_{f_0}(\phi)= h_{\mu}(f_0) + \int \phi \, d\mu,
\quad\mbox{where}\quad
P_{f_0}(\phi)=\sup_{\nu \in \cM}
\left\{ h_{\nu}(f_0) + \int \phi \, d\nu \right\},
\]
and $\cM$ is the set of all $f_0$-invariant probability measures.
The quantity $P_{f_0}(\phi)$ is called the \emph{topological
  pressure} and the identity on the right hand side is a
consequence of the \emph{variational principle}, see
e.g. \cite{Wa82} for definitions of  entropy $h_\mu(f_0)$
and topological pressure $P_{f_0}(\phi)$.

For uniformly expanding maps it turns out that physical
measures are the equilibrium states for the potential
function $\phi (x) = - \log | \det Df_0 (x)|$.  It is a
remarkable fact that \emph{for uniformly hyperbolic and
  uniformly expanding systems these two classes of measures
  (physical and equilibrium states) coincide}, see e.g.
\cite{Bo75,BR75,Pe77,PS82,Ru76}. This relation has been
extended for other types of systems, see e.g.
\cite{Le81,buzzi99,AlOlTah}.

\subsection{The Entropy Formula}
\label{sec:entropy-formula}

Pesin's Entropy Formula \cite{Pe77,LY85,Li98,Li99,QZ2002}
ensures, in particular, that for $C^{1+\alpha}$ maps,
$\alpha>0$, the metric entropy with respect to an invariant
measure $\mu$ \emph{with positive Lyapunov exponents in
  every direction for $\mu$ almost all points} satisfies the
relation
\begin{equation}
  \label{eq:entropyFormula}
  h_\mu(f_0)=\int \log | \det Df_0 (x)| \, d\mu (x)
\end{equation}
\emph{if, and only if, $\mu$ is absolutely continuous with
  respect to the reference measure $m$.}  In general we
integrate the sum of the positive Lyapunov exponents, see
\cite{Li99} for a proof in the $C^2$ setting. In our
setting, the proof that $\mu\ll m$ implies the Entropy
Formula \eqref{eq:entropyFormula} is an exercise using the
bounded distortion provided by the H\"older condition on the
derivative.

We recall that by the Ergodic Theorem any ergodic absolutely
continuous $f_0$-invariant measure $\mu$ is a physical
measure.

Recently \cite{CoYo2004,ArTah,ArTah2} the use of zero-noise
limit measures, obtained through a suitable random
perturbation, to obtain equilibrium states which are then
shown to be $SRB$ or physical measures through the Entropy Formula
has provided many new results (for some earlier work in this
direction see \cite{tsujii92}), enlarging the range of
applicability of smooth ergodic theory.

The relation between equilibrium states and physical
measures provided by the Entropy Formula, among many other
facts, motivates the search for conditions guaranteeing
existence of equilibrium states.
We obtain sufficient
conditions for existence and continuous variation of
equilibrium states in what follows, but first we present
some applications of these results to non-uniformly
expanding systems with criticalities, and to non-uniformly
expanding local diffeomorphisms.



\subsection{Maps with critical sets}\label{s.viana}

Let $f_0: M\to M$ be a $C^2$ map of a compact manifold
$M$ such that its critical set
\[
\cC=\{ x\in M : \det Df_0(x) = 0 \}
\]
is non-degenerate: if $\Delta(x)=\det Df_0(x)$ and
$x\in\cC$, then $\Delta'(x)\neq0$.  Hence
$\cC=\Delta^{-1}(\{0\})$ is a co-dimension one compact
submanifold and thus has zero volume.  We assume further
that $f_0$ is {\em non-flat} near $\cC$: there exist $B>1$
and $\be>0$ for which
 \begin{itemize}
 \item[(S1)]
\hspace{.1cm}$\displaystyle{\frac{1}{B}\dist(x,\cC)^{\be}\leq
\frac{\|Df_0(x)v\|}{\|v\|}\leq B\dist(x,\cC)^{\be}}$;
 \item[(S2)]
\hspace{.1cm}$\displaystyle{\left|\log\|Df_0(x)^{-1}\|-
\log\|Df_0(y)^{-1}\|\:\right|\leq
B\frac{\dist(x,y)}{\dist(x,\cC)^{\be}}}$;
 \item[(S3)]
\hspace{.1cm}$\displaystyle{\left|\log|\det Df_0(x)^{-1}|-
\log|\det Df_0(y)^{-1}|\:\right|\leq
B\frac{\dist(x,y)}{\dist(x,\cC)^{\be}}}$;
 \end{itemize}
 for every $x,y\in M\setminus \cC$ with
$\dist(x,y)<\dist(x,\cC)/2$ and $v\in T_x M$. Given
$\delta>0$ we define the $\delta$-{\em truncated distance}
from $x\in M$ to $\cC$ by
 $$ \dist_\delta(x,\cC)=
\left\{
\begin{array}{ll}
1 & \mbox{if }\dist(x,\cC)\geq \delta,
\\
\dist (x,\cC) & \mbox{otherwise.}
\end{array}
\right.
$$
We observe that if $M$ is one-dimensional (either the interval or
the circle) and $\cC$ is discrete, then (S1)-(S3) amount to
the zeroes of $f_0'$ being non-flat, see~\cite{MS93}.

We assume that $f_0$ is a non-uniformly expanding map, that is
 there is $c>0$ such that
 \begin{equation} \label{liminf1}
\limsup_{n\to +\infty}\frac{1}{n}
\sum_{j=0}^{n-1}\log\|Df_0(f_0^j(x))^{-1}\|\leq -c<0
\end{equation}
for Lebesgue almost every $x\in M$ (recall that we are
taking $\cC$ with zero Lebesgue measure).  Moreover, we
suppose that the orbits of $f_0$ have {\em slow
  approximation to the critical set}, i.e., given small
$\gamma >0$ there is $\delta >0$ such that
 \begin{equation}
 \label{limsup1}
\limsup_{n\to +\infty}\frac{1}{n}
\sum_{j=0}^{n-1}-\log\dist_\delta(f_0^j(x), \cC)\leq \gamma
 \end{equation}
 for Lebesgue almost every $x\in M$. These asymptotic
 conditions are motivated by the following result, which
 non-trivially extends to higher dimensions related results
 from one-dimensional dynamics.

\begin{theorem}[Theorem C in \cite{ABV00}]
   \label{thm:abv-Critical}
   Let $f_0$ be a $C^2$-non-uniformly expanding map whose
   orbits have slow approximation to the critical set. Then
   there is a finite number of ergodic absolutely continuous
   (\emph{SRB}) $f_0$-invariant probability measures
   $\mu_1,\dots ,\mu_p$ whose basins cover a full Lebesgue
   measure subset of $M$, i.e.
   \[
   B(\mu_1)\cup\dots \cup  B(\mu_p) = M, \, m\bmod0.
   \]
   Moreover, every absolutely continuous $f_0$-invariant
   probability measure $\mu$ may be written as a convex
   linear combination of $\mu_1,\dots ,\mu_p$: there are
   non-negative real numbers $w_1,\dots,w_p$ with
   $w_1+\cdots+w_p=1$ for which
   $\mu=w_1\cdot\mu_1+\cdots+w_p\cdot\mu_p$.
\end{theorem}


\subsubsection{Random perturbations and stationary measures}
\label{sec:rand-pert-stat}

Let $\hat f=\{f_t: Y\to Y, t\in X\}$ be a parameterized
family of maps where $X,Y$ are connected compact metric
spaces, which we assume are subsets of finite dimensional
manifolds. We assume that $\hat f: X\times Y\to Y,
(t,x)\mapsto f_t(x)$ is continuous. We consider the
\emph{random iteration of $f$}
\[
f_\omega^n=f_{\omega_n}\circ \dots \circ f_{\omega_1}
\]
for any sequence $\omega=(\omega_1,\omega_2,\dots)$ of
parameters in $X$ and for all $n\ge1$. We let also
$(\theta_\ep)_{\ep>0}$ be a family of \emph{non-atomic}
probability measures on $X$ such that $\supp
(\theta_\ep)\to\{0\}$ when $\ep\to0$.

We set $\Omega=X^\bN$ with the standard infinite product
topology, which makes $\Omega$ a compact metrizable space,
and also take the standard product probability measure
$\theta^\ep=\theta_\ep^\bN$, which makes
$(\Omega,\cB,\theta^\ep)$ a probability space, where $\cB$
is the $\sigma$-algebra of $\Omega$ generated by cylinder
sets.
The following skew-product map is the natural setting for
many definitions connecting random with standard dynamical
systems
\[
   F : \Omega\times Y  \to \Omega\times Y \quad
   (\omega,x) \mapsto (\sigma(\omega),f_{\omega_1}(x))
\]
where $\sigma$ is the left shift on $\Omega$. A probability
measure $\mu^\ep$ on $Y$ is a \emph{stationary measure for
  the random system $(\hat f,\theta_\ep)$} if $\theta^\ep
\times \mu^\ep$ on $\Omega \times Y$ is $F$-invariant.  We
say that $\mu^\ep$ is ergodic if $\theta^{\ep} \times
\mu^\ep$ is $F$-ergodic.

In this setting ($\hat f$ continuous) it is well known that
there always exist an ergodic stationary probability measure
$\mu^\ep$ for all $\ep>0$, see
e.g. \cite{brin-kifer1987}. Moreover every weak$^*$
accumulation point $\mu$ of $(\mu^\ep)_{\ep>0}$ when $\ep\to0$ is
a $f_0$-invariant probability measure, see
e.g. \cite{vdaraujo2001}.

This suggests the notion of stochastic stability: we say
that a map $f_0$ having physical measures (at most countably
many by definition of ergodic basin) is \emph{stochastically
  stable under the perturbation $(\hat f,
  (\theta_\ep)_{\ep>0})$} if \emph{every weak$^*$
  accumulation point $\mu$ of $(\mu^\ep)_{\ep>0}$ when
  $\ep\to0$ is a convex linear combination of the physical
  measures of $f_0$}.


\subsubsection{Non-uniform expansion and slow approximation on
  random orbits}
\label{sec:non-unif-expans}

We study perturbations of endomorphisms by
considering families $(f_t)_{t\in X}$ of maps
where $X$ is a compact connected subset of an Euclidean
space and $f(t,x)=f_t(x)$ is a $C^2$-map satisfying extra
regularity assumptions
\begin{itemize}
\item[(S4)] for all
  $x\in M\setminus\cC$ and every $t,s\in X$ we have
\[
\big| \log|\det Df_t(x)| - \log|\det Df_s(x)| \big| \le
B\cdot d(t,s)^\beta;
\]
\item[(S5)] there exists $\delta_0>0$ such that
\[
\frac{1}{B}\dist(x,\cC)^{\be}\leq
\frac{\|Df_t(x)v\|}{\|v\|}\leq B\dist(x,\cC)^{\be}
\]
for all $x\in M\setminus\cC$ and $d(t,0)<\delta_0$.
\end{itemize}
We further assume that the family of probability measures
$(\theta_\epsilon)_{\epsilon>0}$ is \emph{non-degenerate}:
for all $x\in M$
\begin{itemize}
\item $f(\bullet,x)_* \theta_\epsilon \ll m$ and
\item $\supp f(\bullet,x)_* \theta_\epsilon$ contains a
  neighborhood of $f_0(x)$.
\end{itemize}
This may be implemented, e.g in parallelizable manifolds
(with an additive group structure: tori $\bT^d$ or cylinders
$\bT^{d-k}\times\bR^k$) by considering $ X=\{t\in
\bR^d\colon \|t\|\leq \ep_0\} $ for some $\ep_0>0$,
$\te_\ep$ the normalized Lebesgue measure on the ball of
radius $\ep\leq \ep_0$, and taking $f_t=f+t$; that is,
adding a jump $t$ to the image of $f_0$, which we call
\emph{additive random perturbations}. In general every $C^2$
endomorphism of a compact manifold can be included in a
smooth family of endomorphisms with a non-degenerate family
of probability measures $(\theta_\epsilon)_{\epsilon>0}$ as
above, see e.g. \cite{vdaraujo2000}.

We consider an analog of condition~\eqref{liminf1} for
random orbits. We say that the map $f_0$ is {\em
  non-uniformly expanding for random orbits} if there exists
$c>0$ such that for $\ep>0$ small enough
and for $\te^\ep\times m$ almost every $(\omega,x)\in
\Omega\times M$
 \begin{equation} \label{liminf2}
\limsup_{n\to +\infty}\frac{1}{n}
\sum_{j=0}^{n-1}\log\|Df_0(f_{\omega}^j(x))^{-1}\|\leq -c<0 .
 \end{equation}
 We also consider an analog of condition~(\ref{limsup1}) for
 random orbits; we assume {\em slow approximation of random
   orbits to the critical set}, i.e. given any small $\gamma
 >0$ there is $\delta >0$ such that for $\te^\ep\times m$
 almost every $(\omega,x)\in \Omega\times M$ and small
 $\ep>0$
\begin{equation}
 \label{limsup2}
\limsup_{n\to +\infty}\frac{1}{n}
\sum_{j=0}^{n-1}-\log\dist_\delta(f^j_{\omega}(x), \cC)\leq
\gamma .
 \end{equation}
 Under these conditions we are able to obtain a result on
 the existence of finitely many physical measures for the
 randomly perturbed system.  In the setting of random
 perturbations, a stationary measure $\mu^\ep$ for $(\hat
 f,\theta_\ep)$ is a \emph{physical measure} if its
 \emph{ergodic basin} $B(\mu^\ep)$ has positive Lebesgue
 measure, where
 \begin{align*}
B(\mu^\ep)=\big\{x\in M:
&
\frac1n\sum_{j=0}^{n-1}\varphi(f_\omega^j(x))\to\int\varphi\,
d\mu
\mbox{  when  } n\to\infty \mbox{  for all}
\\
&\mbox{continuous } \varphi: M\to\bR
\mbox{ and }
\theta^\ep\mbox{-almost every } \omega\in\Omega \big\} .
 \end{align*}

\begin{theorem}[Theorem C in~\cite{AA03}]
\label{t.finite2}
Let $f_0: M\to M$ be a $C^2$-non-uniformly expanding map
non-flat near $\cC$, and whose orbits have slow approximation to
$\cC$. If $f_0$ is non-uniformly expanding for random orbits
and random orbits have slow approximation to $\cC$, then
for sufficiently small $\ep>0$ there are physical
 measures $\mu^\ep_1,\dots, \mu^\ep_l$ (with $l$ not
 depending on $\ep$) such that:
 \begin{enumerate}
\item  for each  $x\in M$ and
$\te^\ep$ almost every $\omega\in\Omega$, the average of
Dirac measures $\delta_{f_{\omega}^n(x)}$ converges in the
weak$^*$ topology to some  $\mu^\ep_i$ with $1\le i\le l$;
\item for each $1\le i\le l$ we have
 $$
 \mu^\ep_i=\lim\frac1n\sum_{j=0}^{n-1}\int\big(f_{\omega}^j\big)_*
 \big(m\mid B(\mu^\ep_i)\big)\,d\te^\ep(\omega),
 $$
 where $m\mid B(\mu^\ep_i)$ is the normalization of the
 Lebesgue measure restricted to $B(\mu^\ep_i)$;
\item if $f_0$ is topologically transitive, then $l=1$.
 \end{enumerate}
\end{theorem}

Using Theorem~\ref{t.finite2} together with the general
results from Section~\ref{sec:semic-entr-press} provides a
result on existence of physical measures.

Let $f_0$ be a
non-uniformly expanding $C^2$ map away from a non-flat
critical set $\cC$ and whose orbits have slow approximation
to $\cC$.  We say that $f_0$ is a \emph{stochastically
  robust non-uniformly expanding map} if there exists a
continuous family $\hat f= (f_t)_{t\in X}$ of $C^2$-maps
(with $0\in X$) and a family
$\hat\theta=(\theta_\ep)_{\ep>0}$ of probability measures on
$X$ such that $(\hat f,\hat\theta)$ is non-uniformly
expanding for random orbits and random orbits have slow
approximation to $\cC$.

\begin{theorem}
\label{t.sts}
Let $f_0:M\to M$ be non-uniformly expanding $C^2$ map away
from a non-flat critical set $\cC$ and whose orbits have
slow approximation to $\cC$.

If $f_0$ is a stochastically robust non-uniformly expanding
map, then every weak$^*$ accumulation point $\mu$ of any
family $(\mu^\ep)_{\ep>0}$ of stationary measures given by
Theorem~\ref{t.finite2} is an absolutely continuous
$f_0$-invariant probability measure. In particular $f_0$
admits an absolutely continuous invariant measure.
\end{theorem}

We remark that in the setting of $C^2$ endomorphisms of a
compact manifold, if $\mu$ is an absolutely continuous
$f_0$-invariant measure as in Theorem~\ref{t.sts}, then
$\big|\log|\Delta|\big|$ is $\mu$-integrable,
see~\cite{Li98}. Hence the Entropy Formula holds showing that
$\mu$ is an equilibrium state for $f_0$ with respect to the
potential $-\log|\Delta|$.

In the setting of non-uniformly expanding maps for random
orbits with $m(\cC)=0$ we have that any stationary measure
$\mu^\ep$ from Theorem~\ref{t.finite2} satisfies a similar
Entropy Formula, see e.g.  \cite{BhLi98}
\[
h_{\mu^\ep}=\int \phi_\ep \, d\mu^\ep
\quad\mbox{where}\quad
\phi_\ep(x)=\int \log|\det Df_\omega(x)| \,
d\theta^\ep(\omega).
\]
Putting this together with the abstract results from
Section~\ref{sec:semic-entr-press} enables us to prove the
following weak$^*$ continuity result.

\begin{theorem}
\label{thm.contNUE}
Let $\cW$ be a subset of
non-uniformly expanding $C^2$ maps $f$, with the same
exponent bound $c$, non-flat
near its critical set $\cC(f)$ and whose
orbits have slow approximation to $\cC(f)$.

If we assume that every $f\in\cW$
\begin{enumerate}
\item is a stochastically robust non-uniformly expanding map, and
\item admits a unique ergodic absolutely continuous
$f$-invariant measure $\mu_f$, i.e., $m(M\setminus
B(\mu_f))=0$;
\end{enumerate}
then
\begin{itemize}
\item $\mu_f$ is stochastically stable for every $f\in \cW$, and
\item $\mu_f$ varies continuously with $f\in\cW$
in the weak$^*$ topology.
\end{itemize}
\end{theorem}

Combining Theorems~\ref{t.finite2}, \ref{t.sts} and
\ref{thm.contNUE} we deduce conditions under which
\emph{stochastic stability implies the continuous
  variation of physical measures.}

\begin{theorem}
  \label{thm:existcontNUE}
Let $\cW\subset\diff^2(M)$ be a subset of
non-uniformly expanding $C^2$ maps $f$, with the same
exponent bound $c$, non-flat
near its critical set $\cC(f)$ and whose
orbits have slow approximation to $\cC(f)$.

If $\cW$ is a robustly transitive class (i.e. every
$f\in\cW$ admits a dense orbit) of stochastically robust
non-uniformly expanding maps, then there exists a unique
absolutely continuous ergodic probability measure $\mu_f$
for each $f\in\cW$, satisfying $m(M\setminus B(\mu_f))=0$,
which is stochastically stable and depends continuously on
$f\in\cW$ in the weak$^*$ topology.
\end{theorem}


\subsubsection{Robust non-uniformly expanding maps of the cylinder}
\label{sec:non-trivial-example}

As an application of the previous theorem to the class of
maps on the cylinder $\bS^1\times \bR$ introduced
in~\cite{Vi97}, we obtain with different proofs a version of
results in subsequent works~\cite{Al00} and
\cite{alves-viana2002} where it was shown that such maps
have a unique physical measure which varies continuously
with the map. Here we only provide continuous variation in
the weak$^*$ topology, while the above-mentioned works (much
harder and longer) prove the continuous variation of the
density of the physical measure in the $L^1$ topology.

The class of non-uniformly expanding maps with critical sets
introduced by M. Viana can be described as follows. Let
$a_0\in(1,2)$ be such that the critical point $x=0$ is
preperiodic for the quadratic map $Q(x)=a_0-x^2$. Let
$\bS^1=\bR/\bZ$ and $b:\bS^1\to \bR$ be a Morse function,
for instance, $b(s)=\sin(2\pi s)$. For fixed small
$\alpha>0$, consider the map
 \[
\tilde f:  \bS^1\times\bR
\to  \bS^1\times \bR, \quad
 (s, x) \mapsto \big(\hat g(s),\hat q(s,x)\big)
\]
where $\tilde g$ is the uniformly expanding map of the circle
defined by $\tilde{g}(s)=ds$ (mod $\bZ$) for some $d\ge16$,
and $\tilde q(s,x)=a(s)-x^2$ with $a(s)=a_0+\al b(s)$. It is
not difficult to check that for small enough $\al>0$ there
is an interval $I\subset (-2,2)$ such that $\tilde
f(\bS^1\times I)\subset\inte( S^1\times I)$. Hence every map
$f$ $C^0$-close to $\tilde f$ has $\bS^1\times I$ as a forward
invariant region.  We consider these maps $f$ close to $\tilde
f$ restricted to $\bS^1\times I$.  By the expression of
$\tilde f$ it is not difficult to verify that $\tilde f$, and
any map $f$ $C^2$-close to it, is non-flat near the
critical set.

\begin{theorem}[Theorem A in \cite{Vi97} $\&$ Theorem C in
  \cite{alves-viana2002} $\&$ Theorem E in \cite{AA03}]
\label{thm.AA}
If $f$ is sufficiently close to $\tilde f$ in the $C^3$
topology then $f$ is topologically mixing, non-uniformly
expanding and its orbits have slow approximation to the
critical set.  Moreover if the noise level $\ep$ of an
additive random perturbation $(\hat f,\theta_\ep)$ of $f$ is
sufficiently small, then $f$ is non-uniformly expanding for
random orbits and random orbits have slow approximation to
the critical set.
 \end{theorem}

 As an immediate consequence of Theorems~\ref{thm.AA}
 and~\ref{thm:existcontNUE} we have that for an open subset
 $\cU$ of $C^3$ maps near $\hat f$ \emph{there is a unique
   physical measure $\mu_f$ for $f\in\cU$ which varies
   continuously with $f\in\cU$ in the weak$^*$ topology.}


 \subsubsection{Quadratic maps}
 \label{sec:quadratic-maps}

For the family $f_a(x)=a-x^2$, $a\in[0,2]$ it is known
\cite{BC85,BC91} that there exists a positive measure subset
$B\subset[0,2]$ such that for $a\in B$ the map $f_a$ is
transitive and non-uniformly expanding with the a uniform
exponent bound $c>0$. In a recent work \cite{freitas} the
orbits of these maps where shown to have slow approximation
to the non-flat critical set $\cC=\{0\}$.

As a consequence of these results and of
Theorem~\ref{thm:existcontNUE}, we have that each map $f_a$,
with $a\in B$, admits an absolutely continuous invariant
probability measure $\mu_a$, the map $a\in B\mapsto\mu_a$ is
weak$^*$ continuous, $B(\mu_a)=[-1,1], \, m \bmod 0$ and
each $\mu_a$ is stochastically stable, $a\in B$.


\subsection{Local diffeomorphisms}
\label{sec:local-diff}
Let $f_0:M\to M$ be a $C^2$ local diffeomorphism of the
manifold $M$ and assume that $f_0$ satisfies
condition~\eqref{liminf1} for Lebesgue almost every $x\in
M$. We are in the setting of maps ``with empty critical set
$\cC=\emptyset$'' so
Theorems~\ref{thm:abv-Critical},~\ref{t.finite2},~\ref{t.sts},
\ref{thm.contNUE} and~\ref{thm:existcontNUE} also hold
since~\eqref{limsup1} and~\eqref{limsup2} are vacuous.

In~\cite{AA03} sufficient conditions and
necessary conditions were obtained for stochastic stability
of non-uniformly expanding local diffeomorphisms. Using
results from Subsection~\ref{sec:rand-pert-upper} on
zero-noise limits of random equilibrium states we obtain a
necessary and sufficient condition for stochastic stability
in this setting.

\begin{theorem}
\label{thm.stc}
Let $f_0: M\to M$ be a non-uniformly
expanding $C^2$ local diffeomorphism. Then $f_0$ is
stochastically stable if, and only if, $f_0$ is non-uniformly
expanding for random orbits.
\end{theorem}

\subsubsection{Equilibrium states for potentials of low
  variation}
\label{sec:equil-stat-potent}

We consider the following class $\cU$ of $C^2$ local
diffeomorphisms $f:M\to M$  which may be seen as small deformations
of uniformly expanding maps. We assume that for positive
constants $\de_0, \beta,\de_1,\sigma_1$ and integers $p,q$
there exists a covering $B_1,\dots,B_{p+q}$ of $M$ such that
$f\mid B_i$ is injective for all $i=1,\dots,p+q$ and
\begin{enumerate}
\item $f$ expands uniformly at $x\in B_1\cup\dots\cup B_p$:
 $ \| Df(x)^{-1} \|\le (1+\de_1)^{-1}$;
\item $f$ never contracts too much: $\| Df(x)^{-1} \|\le
  1+\de_0$ for all $x\in M$;
\item $f$ is volume expanding: $|\det Df(x)|\ge\sigma_1$ for
  all $x\in M$ with $\sigma_1>p$;
\item there exists a set $W$ such that
  \begin{enumerate}
  \item $V=\{ x\in M : \|
  Df(x)^{-1} \|\ge (1+\de_1)^{-1} \}\subset W \subset
  B_{p+1}\cup\dots\cup B_{p+q}$;
\item $\inf \log \|Df\mid M\setminus W\| > \sup \log
  \|Df\mid V\|$; and
  \end{enumerate}
\item $\sup \log \|Df\mid V\| - \inf \log \|Df\mid V\| < \beta$.
\end{enumerate}
We observe that $\cU$ contains an open set of $C^2$ local
diffeomorphisms on tori $\bT^n, n\ge2$, see
e.g.~\cite{ABV00,AA03}.

Given a continuous function $\phi:M\to\bR$ and $\rho>0$ we say that
\emph{$\phi$ has $\rho$-low variation} if
\[
\sup \phi \le P_{f}(\phi) - \rho\cdot\htop(f),
\]
where $\htop(f)$ is the \emph{topological entropy} of $f$
which coincides (through the variational principle) with the
pressure $P_f(0)$ for any constant potential.

\begin{theorem}[\cite{krerley2003}]
  \label{thm:krerley}
  For $\de_0$ and $\beta$ small enough there exists
  $\rho_0>0$ such that, for all $f\in\cU$ and
  $0<\rho<\rho_0$, every $\phi:M\to\bR$ of $\rho$-low
  variation admits an ergodic equilibrium state $\mu_\phi$.
  Moreover $\mu_\phi(\log\| (Df)^{-1}\| ) \le c =
  c(\de_1,\sigma_1,p,q) <0$, that is, every Lyapunov
  exponent of $\mu_\phi$ is positive.
\end{theorem}

We note that \emph{ the notion of low variation potential
  includes the constant potentials}.  Hence for this
$C^2$-open class $\cU$ of maps there are measures of maximal
entropy, which are equilibrium states for the potential
$\phi\equiv0$. We may apply to these maps the abstract
Theorem~\ref{thm:semicontpressure} from
Section~\ref{sec:semic-entr-press} to deduce the following.

\begin{theorem}
  \label{thm:semiconthtop}
  When restricted to maps in $\cU$, topological entropy
  $\htop: \cU\to\bR, f\mapsto \htop(f)$ is
  an upper semicontinuous function.
\end{theorem}

In Section~\ref{sec:semic-entr-press} we present the
abstract results used to prove Theorems~\ref{t.sts},
\ref{thm.contNUE}, \ref{thm.AA} and \ref{thm.stc}. In
Section~\ref{sec:proof-semic-meas} we prove the abstract
results. Finally in Section~\ref{sec:stoch-stab-crit} we
show how to derive the above-mentioned theorems from the
results in Section~\ref{sec:semic-entr-press}.


\section{Semicontinuity of pressure,  entropy and
  equilibrium states}
\label{sec:semic-entr-press}

Now we state the main technical results.
In the following statements $X,Y$ denote compact metric spaces.

Given a map $f:Y\to Y$ and a Borel probability measure $\mu$
we say that a $\mu\bmod 0$ partition $\xi$ of $Y$ is a
\emph{generating partition} if
\[
\bigvee_{i=0}^{+\infty} (f^i)^{-1}\xi = \cA, \quad \mu\bmod0,
\]
where $\cA$ is the Borel $\sigma$-algebra of $Y$. We denote
by $\partial\xi$ the set of topological boundaries of all
elements of $\xi$.

\begin{theorem}[Upper semicontinuity of topological pressure]
\label{thm:semicontpressure}
Let $f:X\times Y\to Y$ define a  family of
continuous maps $f_t:Y\to Y, y\in Y\mapsto f_t(y)=f(t,y)$
and $(\phi_t)_{t\in X}$ a family of continuous functions
(potentials) $\phi_t: Y\to\bR$ satisfying the following
conditions.
\begin{enumerate}
\item $f_t$ admits some equilibrium state for $\phi_t$, i.e.
  there exists $\mu_t\in\cP_{f_t}(Y)$ such that
  $P_{f_t}(\phi_t)=h_{\mu_t}(f_t)+\int \phi_t \, d\mu_t$ for
  all $t\in X$.
\item Given a weak$^*$ accumulation point $\mu_0$ of $\mu_t$
  when $t\to 0\in X$, let $t_k\to0$ when $k\to\infty$ be
  such that $\mu_k=\mu_{t_k}\to\mu_0$. We write
  $f_k=f_{t_k}, \, \phi_k=\phi_{t_k}$ and assume also that
  \begin{enumerate}
  \item $f_k(y)\to f_0(y)$  when $k\to\infty$ for all $y\in Y$.
  \item there exists a finite $\mu_k$-modulo zero partition
  $\xi$ of $Y$ which is generating for $(Y,f_k,\mu_k)$, $k\ge1$, and
  $\mu_0(\partial \xi)=0$.
  \item $\limsup_{k\to\infty} \int \phi_k\, d\mu_k \le \int \phi_0
  \, d\mu_0$.
  \end{enumerate}
\end{enumerate}
Then $\limsup_{k\to\infty} P_{f_k}(\phi_k)\le P_{f_0}(\phi_0)$.
\end{theorem}

Theorem~\ref{thm:semicontpressure} is a simple consequence
of the next result.

\begin{theorem}[Upper semicontinuity of measure-theoretic entropy]
  \label{thm.semicontentropy}
  Let $f_t:Y\to Y$ be a family of continuous maps as above and
  $\mu_t$ a family of $f_t$-invariant probability measures
  for $t\in X$. Given a weak$^*$ accumulation point $\mu_0$
  of $\mu_t$ when $t\to 0\in X$, we let $t_k\to0$ when
  $k\to\infty$ be such that $\mu_k=\mu_{t_k}\to\mu_0$ and
  write $f_k=f_{t_k}$.

  If there exists a finite $\mu_k$-modulo zero partition
  $\xi$ of $Y$ which is generating for $(Y,f_{t_k},\mu_k)$,
  $k\ge1$, and $\mu_0(\partial \xi)=0$, then $
  \limsup_{k\to\infty} h_{\mu_k}(f_k) \le h_{\mu_0}(f,\xi).$
\end{theorem}

>From this we easily deduce the following.

\begin{theorem}[Continuity of equilibrium states]
  \label{thm:contpressure}
  Let $f_t:Y\to Y$ be a family of continuous maps and
  $\phi_t: Y\to\bR$ a family of continuous functions
  (potentials) satisfying conditions 1 and 2 on
  Theorem~\ref{thm:semicontpressure}, for $t\in X$.

  If $P_{f_k}(\phi_k)\to P_{f_0}(\phi_0)$ for a sequence
  $t_k\to0\in X$, then every weak$^*$ accumulation point
  $\mu$ of $(\mu_k)_{k\ge1}$ when $k\to\infty$ is a
  equilibrium state for $f_0$ and the potential $\phi_0$.
\end{theorem}

\subsection{Upper semicontinuity of
  random measure-theoretic entropy}
\label{sec:rand-pert-upper}

We need the notion of metric entropy  for random
dynamical systems which may be defined as follows.

\begin{theorem}[Theorem 1.3 in \cite{Ki86}] \label{thm.metr-entr-rand}
For any finite measurable partition $\xi$ of $Y$ the limit
 $$
 h_{\mu^\ep}((\hat f,\theta_\ep), \xi) =
 \lim_{n \rightarrow \infty}
 \frac{1}{n} \int H_{\mu^\ep} \left(
 \bigvee_{k=0}^{n-1}(f^k_{\omega})^{-1} \xi \right) d
 \theta^{\ep} (\omega)
 $$
 exists. This limit is called the entropy of the random
 dynamical system with respect to $\xi$ and to $\mu^\ep$.
\end{theorem}

As in the deterministic case the above limit can be replaced
by the infimum. The \emph{metric entropy of the random
  dynamical system $(\hat f,\theta_\ep)$} is given by
$h_{\mu^\ep} (\hat f, \theta_\ep)= \sup h_{\mu^\ep} ((\hat
f,\theta_\ep), \xi)$, where the supremum is taken over all
measurable partitions.

Kolmogorov-Sinai's result about generating partitions is
also available for random maps.  We say that $\xi$ is a $\mu^\ep$-random
generating partition if $\xi$ is a finite
partition of $Y$ such that
$$
\bigvee_{k=0}^{+\infty} (f_{\omega}^k)^{-1} \xi = \cA, \,
\mu^\ep\bmod0
\quad \text{for} \quad \theta^{\ep}-\text{almost all }
\omega \in \Omega.
$$

\begin{theorem}[Corollary 1.2 in \cite{Ki86}]
  \label{thm.KSrandom}
  If $\xi$ is a $\theta_\ep$-random generating partition,
  then we have $h_{\theta^\ep\times\mu^\ep}(\hat  f,\theta_\ep)=
  h_{\mu^\ep}((\hat f,\theta_\ep),\xi)$.
\end{theorem}

Now we can state the following upper-semicontinuity property.

\begin{theorem}[Upper semicontinuity of random measure-theoretic entropy]
  \label{thm.semicontrandomentropy}
  Let $\mu$ be the weak$^*$ limit of $(\mu^{\ep_k})_{k\ge1}$
  when $k\to\infty$ for a sequence $\ep_k\to0$.  Let us
  assume that there exists a finite partition $\xi$ of $Y$
  which is $\theta_{\ep_k}$-generating for random orbits,
  for every $k\ge1$, and such that $\mu(\partial\xi)=0$.
  Then
\[
\limsup_{k\to\infty} h_{\mu^{\ep_k}}((\hat f, \te_{\ep_k}),\xi)
\le h_{\mu}(f,\xi).
\]
\end{theorem}

As a consequence of this we deduce a result which provides a
way to obtain equilibrium states as zero-noise limits.

\begin{theorem}[Continuity of random equilibrium states]
  \label{thm:zero-noise-equilibrium}
  Let $\mu$ be the weak$^*$ limit of
  $(\mu_k=\mu^{\ep_k})_{k\ge1}$ when $k\to\infty$ for a
  sequence $\ep_k\to0$ when $k\to\infty$. Let us assume that
  there exists a finite partition $\xi$ of $Y$ $\mu_k\bmod0$
  which is $\theta_{\ep_k}$-generating for random orbits for
  all $k\ge1$.

  Moreover we suppose that $h_{\mu_k}(\hat
  f,\theta_{\ep_k})=\int \phi_k \,d\mu_k$ for all $k\ge1$
  where $\phi_k:Y\to\bR$ is a sequence of functions such
  that $\phi_k\to\phi_0$ pointwisely when $k\to\infty$ and
  $P_{f_0}(-\phi_0)\le0$.  Then $\mu$ is an equilibrium
  state for $-\phi_0$, that is $h_\mu(f_0)=\int\phi_0\,d\mu$.
\end{theorem}



\section{Proof of semicontinuity of measure-theoretic
  entropy and equilibrium states}
\label{sec:proof-semic-meas}

\subsection{The random setting}
\label{sec:random-setting}

Here we prove Theorem~\ref{thm.semicontrandomentropy} and
Theorem~\ref{thm:zero-noise-equilibrium}.

Let $(\hat f,(\theta_\ep)_{\ep>0})$ be a random perturbation
of $f_0:Y\to Y$, $\mu^0$ be the weak$^*$ limit of
$(\mu^{\ep_k})_{k\ge1}$ when $k\to\infty$ for a sequence
$\ep_k\to0$ and let $\xi$ be a finite
$\theta_{\ep_k}$-generating partition for random orbits, for
all $k\ge1$, as in the statement of
Theorem~\ref{thm.semicontrandomentropy}, that is
$\mu(\partial\xi)=0$.

We first construct a sequence of partitions of
$\Omega$ according to the following result. For a
partition $\cP$ and $y\in \Omega$ we denote
by $\cP(y)$ the element (atom) of $\cP$ containing $y$.
We set $\omega_0=(0,0,0,\dots)\in\Omega$ in what follows.

\begin{lemma}
\label{lem.seqnpartitions} There exists an increasing
sequence of measurable partitions $(\cB_n)_{n\ge1}$ of
$\Omega$ such that
 \begin{enumerate}
 \item $\omega_0\in\inte(\cB_n(\omega_0))$ for all $n\ge1$;
 \item $\cB_n \nearrow \cB$, $\theta^{\ep_k} \bmod 0$ for all
   $k\ge1$ when $n\to\infty$;
 \item $\lim_{n \to \infty} H_{\rho} (\xi \mid \cB_n) =
   H_{\rho} (\xi \mid \cB)$ for every measurable finite
   partition $\xi$ of $\Omega$ and any $F$-invariant probability measure
   $\rho$.
 \end{enumerate}
\end{lemma}

\begin{proof}
For the first two items we let $\cE_n$ be a finite
$\theta_{\ep_k}\bmod 0$ partition of $X$ such that
$0\in\inte(\cE_n(t_0))$ with $\diam(\cE_n)\to0$ when
$n\to\infty$. Example: take a cover $(B(t,1/n))_{t\in X}$ of
$X$ by $1/n$-balls and take a subcover $U_1,\dots,U_k$ of
$X\setminus B(t_0,2/n)$ together with $U_0=B(t_0,3/n)$; then
let $\cE_n=\{U_0,M\setminus
U_0\}\vee\dots\vee\{U_k,M\setminus U_k\}$.

We observe that we may assume that the boundary of these
balls has null $\theta_{\ep_k}$-measure for all $k\ge1$, since
$(\theta_{\ep_k})_{k\ge1}$ is a denumerable family of
non-atomic probability measures on $X$ and $X$ may be taken
as a subset of some Euclidean space.  Now we set
\[
\cB_n=\cE_n\times\stackrel{n}{\dots}\times\cE_n\times\Omega\quad
\mbox{for all  }n\ge1.
\]
Then since $\diam \cE_n\le 2/n$ for all $n\ge1$ we have that
$\diam \cB_n\le 2/n$ also and so tends to zero when
$n\to\infty$.  Clearly $\cB_n$ is an increasing sequence of
partitions. Hence $\vee_{n\ge1} \cB_n$ generates the
$\sigma$-algebra $\cB$, $\theta^{\ep_k} \bmod0$ (see e.g.
Lemma 3 of Chapter 2 in \cite{Bi65}) for all $k\ge1$. This
proves items (1) and (2).

Item 3 of the statement of the lemma is Theorem 12.1 in~\cite{Bi65}.
\end{proof}

Now we use some known properties of conditional entropy to
derive the right inequalities. First we recall that
\begin{eqnarray*}
  h_{\mu^{\ep_k}}(\hat f,\theta_{\ep_k})
  &=&
  h_{\mu^{\ep_k}}((\hat f,\theta_{\ep_k}),\xi)
  =
  h_{\theta^{{\ep_k}} \times \mu^{\ep_k}}^{ \cB \times Y} (F,
  \Omega\times\xi)
  \\
  &=&
  \inf\frac1n H_{\theta^{\ep_k}\times\mu^{\ep_k}}\left(
  \bigvee_{j=0}^{n-1} (F^j)^{-1}(\Omega\times\xi) \mid
  \cB\times Y \right)
\end{eqnarray*}
where the first equality comes from the Kolmogorov-Sinai
Theorem~\ref{thm.KSrandom} and the assumption that $\xi$ be
generating, while the second one can be found in Theorem 1.4
of Chapter II in~\cite{Ki86}, with $\Omega\times\xi=\{
\Omega\times A: A\in\xi\}$.  Here $ h_{\theta^{{\ep_k}}
  \times \mu^{\ep_k}}^{ \cB \times Y} (F, \Omega\times\xi)$
is the \emph{conditional entropy of $\theta^{{\ep_k}} \times
  \mu^{\ep_k}$ with respect to the $\sigma$-algebra $\cB
  \times Y$ on the partition $\Omega\times\xi$}, whose
definition is translated in the second line of the above
formula and whose basic properties can be found
in~\cite{Bi65,Ki86}.

The last expression shows that for arbitrary fixed $N\ge1$ and
for any $l\ge1$
\begin{eqnarray*}
 h_{\mu^{\ep_k}}(\hat f,\theta_{\ep_k})
 &\le&
 \frac1N H_{\theta^{\ep_k}\times\mu^{\ep_k}}\left(
  \bigvee_{j=0}^{N-1} (F^j)^{-1}(\Omega\times\xi) \mid
  \cB\times Y \right)
 \\
 &\le&
 \frac1N H_{\theta^{\ep_k}\times\mu^{\ep_k}}\left(
  \bigvee_{j=0}^{N-1} (F^j)^{-1}(\Omega\times\xi) \mid
  \cB_l\times Y \right)
\end{eqnarray*}
because $\cB_l\times Y\subset\cB\times Y$. Now we fix $N$ and $l$,
let $k\to \infty$ and note that since
$\mu^0(\partial\xi)=0=\de_{\omega_0}(\partial \cB_l)$ it must be
that
\[
(\de_{\omega_0}\times\mu^0)(\partial(B_i\times\xi_j))=0\quad
\mbox{for all  } B_i\in\cB_l \mbox{  and  } \xi_j\in\xi,
\]
where $\de_{\omega_0}$ is the Dirac mass concentrated at
$\omega_0\in\Omega$.  Thus we get by weak$^*$ convergence of
$\theta^{\ep_k}\times\mu^{\ep_k}$ to
$\delta_{\omega_0}\times\mu^0$ when $k\to\infty$
\begin{equation}
  \label{eq:5}
  \limsup_{k\to\infty} h_{\mu^{\ep_k}}(\hat
  f,\theta_{\ep_k})
  \le
  \frac1N H_{\de_{\omega_0}\times\mu^0}\left(
  \bigvee_{j=0}^{N-1} (F^j)^{-1}(\Omega\times\xi) \mid
  \cB_l\times M \right).
\end{equation}
Here it is easy to see that the conditional entropy on the
right hand side of~\eqref{eq:5} (involving only finite
partitions) equals
\begin{equation}
  \label{eq:6}
\frac1N H_{\mu^0}\big(\bigvee_{j=0}^{N-1} f^{-j}\xi \big)
=\frac1{N}\sum_i \mu^0(P_i)\log\mu^0(P_i),
\end{equation}
with $P_i=\xi_{i_0}\cap f^{-1}\xi_{i_1}\cap\dots\cap
f^{-(N-1)} \xi_{i_{N-1}}$ ranging over all possible
sequences $\xi_{i_0},\dots,\xi_{i_{N-1}}$ of elements of
$\xi$.

Finally, since $N$ was an arbitrary integer,
Theorem~\ref{thm.semicontrandomentropy} follows
from~\eqref{eq:5} and~\eqref{eq:6}.

\medskip

Now to prove Theorem~\ref{thm:zero-noise-equilibrium} we
assume in addition that for each $\mu_k$ there exists a
continuous potential $\phi_k:Y\to\bR$ such that $h_{\mu_k}(\hat
f,\theta_{\ep_k})=\int \phi_k \,d\mu_k$, for $k\ge1$.
Moreover $\phi_k\to\phi_0$ pointwisely to a continuous
potential $\phi_0$ when $k\to\infty$ and
$P_{f_0}(-\phi_0)\le0$.
Then by the previous arguments
\begin{equation*}
\int\phi_0\,d\mu_0
=
\limsup_{k\to\infty} h_{\mu^{\ep_k}}(\hat
f,\theta_{\ep_k})
\le
h_{\mu_0}(f_0,\xi)
\le
h_{\mu_0}(f_0)
\le
\int\phi_0\,d\mu_0
\end{equation*}
concluding the proof of
Theorem~\ref{thm:zero-noise-equilibrium}.

\subsection{The deterministic setting}
\label{sec:determ-sett}

Here we prove
Theorems~\ref{thm:semicontpressure},~\ref{thm.semicontentropy}
and~\ref{thm:contpressure}.

Let $f_k:Y\to Y$ be a sequence of continuous maps such that
$\mu_k$ is $f_k$-invariant for all $k\ge1$, $f_0:Y\to Y$
is continuous with $f_k\to f_0$ pointwisely and
$\mu_k\to\mu_0$ in the weak$^*$ topology when $k\to\infty$.
Let $\xi$ be a finite $\mu_k$-modulo zero partition $\xi$ of
$Y$ which is generating for $(Y,f_{t_k},\mu_k)$, $k\ge1$,
and $\mu_0(\partial \xi)=0$.

Following the same reasoning as
in Subsection~\ref{sec:random-setting} we have for any
given fixed $N\ge1$ that
\[
h_{\mu_k}(f_k)=h_{\mu_k}(f_k,\xi)=
\inf_{n\ge1}\frac1n H_{\mu_k}\big(
\vee_{j=0}^{n-1} (f_k^j)^{-1} \xi \big)
\le\frac1{N}  H_{\mu_k} (\xi_k^N ),
\]
since $\xi$ is generating, where $\xi_k^N=\vee_{j=0}^{N-1}
(f_k^j)^{-1} \xi$. But $\mu_0(\partial\xi)=0$ so for any
given $N\ge1$ we have $\mu_0(\partial\xi_0^N)=0$ also
because $\mu_0$ is $f_0$-invariant. Moreover the weak$^*$
convergence and $f_k$-invariance ensures that
$(f_k^i)_*\mu_k=\mu_k\to\mu_0$ for all $i\ge0$, hence
$\mu_k(\xi_k^N(z))\to\mu_0(\xi_0^N(z))$ when $k\to\infty$
for $\mu_0$-almost every $z\in Y$. In particular we get for
arbitrary $N\ge1$
\[
\limsup_{k\to+\infty} h_{\mu_k}(f_k)
\le
\frac1{N} H_{\mu_0}(\xi_0^N)
\quad\mbox{and so}\quad
\limsup_{k\to+\infty} h_{\mu_k}(f_k)
\le h_{\mu_0}(f_0,\xi)
\]
concluding the proof of Theorem~\ref{thm.semicontentropy}.

\medskip

To prove Theorem~\ref{thm:semicontpressure} we assume in
addition that for each $k\ge1$ there exists a potential
$\phi_k$ and a probability measure $\mu_k$
such that $P_{f_k}(\phi_k)= h_{\mu_k}(f_k) + \int
\phi_k\, d\mu_k$. If we assume also condition (2b) from the
statement of  Theorem~\ref{thm:semicontpressure},
then the result follows  using
Theorem~\ref{thm.semicontentropy} since
\[
\limsup_{k\to\infty} P_{f_k}(\phi_k)
\le
\limsup_{k\to\infty}  h_{\mu_k}(f_k) +
\limsup_{k\to\infty}\int \phi_k\, d\mu_k
\le h_{\mu_0}(f_0) + \int \phi_0\, d\mu_0.
\]

Moreover if we assume that $P_{f_k}(\phi_k)\to
P_{f_0}(\phi_0)$ when $k\to\infty$,
then the same argument above gives
\[
P_{f_0}(\phi_0) \le
\limsup_{k\to\infty}  h_{\mu_k}(f_k) +
\limsup_{k\to\infty}\int \phi_k\, d\mu_k
\le h_{\mu_0}(f_0) + \int \phi_0\, d\mu_0
\le P_{f_0}(\phi_0),
\]
showing that $\mu_0$ is an
equilibrium state for $f_0$ with respect to the potential
$\phi_0$, thus proving Theorem~\ref{thm:contpressure}.


\section{Statistical stability  for non-uniformly
  expanding maps}
\label{sec:stoch-stab-crit}

Here we prove the results in Subsections~\ref{s.viana}
and~\ref{sec:local-diff}.


\subsection{Maps with critical sets}
\label{sec:maps-with-critical}

Here we prove Theorem~\ref{t.sts}.

Let $f_0:M\to M$ be a $C^2$ non-uniformly expanding map away
from a non-flat critical set $\cC$ whose orbits have slow
approximation to $\cC$.  Let also $\hat f= (f_t)_{t\in X}$
be a continuous family in $C^2(M,M)$ and
$\hat\theta=(\theta_\ep)_{\ep>0}$ be a family of probability
measures on $X$ such that $(\hat f,\hat\theta)$ is
non-uniformly expanding for random orbits and random orbits
have slow approximation to $\cC$.

According to Theorem~\ref{t.finite2}, for every small
$\ep>0$ there exists an absolutely continuous stationary
probability measure $\mu^\ep$. Since every $f_t$ is a $C^2$
endomorphism, the random version of the Entropy Formula
ensures that (see e.g. \cite{Li99}) $\mu^\ep$ is an
equilibrium state for $\phi_0=-\log|\det Df_0|$:
\[
h_{\mu^\ep} = \int \int \log| \det Df_t (x) | \, d\theta_\ep(t)
\, d\mu^\ep(x)
\ge c\cdot \dim(M),
\]
since every Lyapunov exponent of the random system is
bounded away from zero by a uniform constant $c>0$
and the sum of all Lyapunov exponents is given by
the above integral.

\begin{lemma}\label{le:unifjacobian}
  Given families $(\hat f, (\theta_\epsilon)_{\epsilon>0})$,
  if random orbits have slow approximation to $\cC$, then
  $\log\big|\det Df_\omega\big|$ is uniformly integrable with
  respect to the family $\theta^\epsilon\times\mu^\epsilon$,
  i.e. given $\gamma>0$ there exists $\delta>0$ such that
  for all small enough $\epsilon>0$ we have
\[
\int_{B(\cC,\delta)} d\mu^\epsilon(x)
\int d\theta^\epsilon(\omega)
-\log\big|\det Df_\omega\big| \le \gamma
\]
and so given a limit point $\mu^0=\lim_{k\to\infty}
\mu^{\epsilon_k}$ in the weak$^*$ topology we have
\[
\int d\mu^{\epsilon_k}(x) \int d\theta^{\epsilon_k}(\omega)
\log\big|\det Df_\omega\big| \to \int\log|\det Df_0| \,
d\mu^0.
\]
\end{lemma}

Here we write $B(A,\delta)=\cup_{x\in A}
B(x,\delta)$ for the $\delta$-neighborhood of a subset $A$.

\begin{proof}
  Since for each small $\epsilon>0$ we have that $\log|\det
  Df_\omega|$ is
  $\theta^\epsilon\times\mu^\epsilon$-integrable, then the
  Ergodic Theorem ensures that
\begin{align}
  \int_{B(\cC,\delta)} d\mu^\epsilon(x) &
\int d\theta^\epsilon(\omega)
-\log\big|\det Df_\omega\big|= \nonumber
\\
&
\lim_{n\to\infty} -\frac1n\sum_{j=0}^{n-1}
\chi_{B(\cC,\delta)}\big( f^j_\omega(x)\big)\cdot
\log\big| \det Df_{\omega_{j+1}}\big( f^j_\omega(x)\big)
\big|
\label{eq:ET}
\end{align}
for $\theta^\epsilon\times\mu^\epsilon$-a.e. $(\omega,x)$.
But the non-degeneracy condition (S4)  ensures that
if we assume $\epsilon>0$ is so small that each
$\omega=(\omega_1,\omega_2, \dots)\in\supp\theta^\epsilon$
satisfies $d(\omega_i,0)<\delta$ for all $i\ge1$, then
(\ref{eq:ET}) is bounded by
\[
\lim_{n\to\infty} -\frac1n\sum_{j=0}^{n-1}
\chi_{B(\cC,\delta)}\big( f^j_\omega(x)\big)\cdot
\log\big| \det Df\big( f^j_\omega(x)\big)
\big| + B\cdot\delta^\beta.
\]
By the non-flatness condition (S3) on $\cC$ there exists a
constant $K>1$ such that
\[
\frac1K\cdot \dist\big( x, \cC \big) \le
\big| \det Df(x) \big| \le
K\cdot \dist\big( x, \cC \big)
\quad\mbox{if}\quad
\dist\big( x, \cC \big)<\delta,
\]
and so we may bound (\ref{eq:ET}) by
\begin{align*}
 & \limsup_{n\to\infty} \frac1n\sum_{j=0}^{n-1}
-\log\big[ K\cdot \dist_\delta\big(
f^j_\omega(x),\cC\big)
\big] + B\cdot\delta^\beta
\\
\le &
\limsup_{n\to\infty} \frac1n\sum_{j=0}^{n-1}
-\log \dist_\delta\big(
f^j_\omega(x),\cC\big) + B\cdot\delta^\beta.
\end{align*}
The assumption of slow approximation to $\cC$ for random
orbits guarantees that the above limit can be made smaller
than $\gamma/2$ for any given $\gamma>0$ by choosing $\delta$ and
$\epsilon$ small enough, and this in turn ensures that the
last expression can be made smaller that $\gamma$. This
concludes the proof.
\end{proof}

According to the non-degeneracy condition (S5) we have that
\[
\log\big\| Df_t(x)^{-1} \big\| \le
\log B -\beta\cdot \log\dist(x,\cC)
\]
therefore using this bound together with the slow
approximation to $\cC$ for random orbits we see that
$\log\big\| Df_t(x)^{-1} \big\|$ is also uniformly
integrable with respect to the family
$\theta^\epsilon\times\mu^\epsilon$.

Now we choose a stationary measure $\mu^{\ep_k}$ for a
sequence $\ep_k\to0$ and take any weak$^*$ accumulation
point $\mu_0$ of $( \mu^{\ep_k} )_k$ when $k\to\infty$.

If we assume that a uniform random generating partition
exists, then by Theorem~\ref{t.sts} and by
Lemma~\ref{le:unifjacobian} we get that $\mu_0$ is
$f_0$-invariant and satisfies
\begin{equation}
  \label{eq:eqstate>0}
h_{\mu_0}(f_0)=\int\phi_0\,d\mu_0\ge c\cdot \dim(M)>0
\quad\mbox{and}\quad
\mu_0\Big(\log\big\| Df_t(x)^{-1} \big\|\Big)\le -c <0.
\end{equation}
But the characterization of measures satisfying the Entropy
Formula for endomorphisms, see e.g. \cite{Li98}, ensures
that a $f_0$-invariant probability measure $\mu_0$
satisfying (\ref{eq:eqstate>0}) is absolutely continuous.

This finishes the proof of
Theorem~\ref{t.sts} except for the
existence of a uniform random generating partition, which is
the content of the following subsection.

\subsubsection{Uniform generating partitions for
  equilibrium states}
\label{sec:unif-rand-gener}

To build a uniform random generating partition for
equilibrium measures we make use of the following notion:
given $0<\al<1$ and $\delta>0$, we say that $n\in\bN$ is a
\emph{$(\al,\delta)$-hyperbolic time for $(\omega,x)\in
  \Omega\times M$} if
\[
\prod_{j=n-k}^{n-1}
\|Df_{\omega_{j+1}}(f^j_{\omega}(x))^{-1}\|\leq\al^k
\quad\mbox{and}\quad
\dist_\delta(f_{\omega}^{n-k}(x),\cC)\geq \al^{bk}
\]
 for every $1\leq k\leq n$, where $\Omega=X^\bN$ was
 defined in Subsection~\ref{sec:rand-pert-stat}. The
 following results ensures the existence of hyperbolic times
 in our setting.

\begin{proposition}[Proposition 2.3 in \cite{AA03}]
  \label{pr:hyptimes}
If $(\hat f, \theta^\ep)$ is non-uniformly expanding for
random orbits and random orbits have slow approximation to
the critical set $\cC$, then there are $\de>0$ and
$\alpha\in(0,1)$ such that $\theta^\ep\times m$-almost
every $(\omega,x)\in\Omega\times M$ has infinitely many
$(\alpha,\de)$-hyperbolic times.
\end{proposition}

\begin{remark}
  \label{rmk:nocritical}
When $\cC=\emptyset$ the second condition on the definition
of hyperbolic time is vacuous and in this case we just write
$\de$-hyperbolic time. Moreover setting
$\omega_t=(t,t,t,\dots)$ then a hyperbolic time for
$x$ with respect to a map $f_t$ is just the same as a
hyperbolic time for $(\omega_t,x)$, $t\in X$.
\end{remark}

Now we state the main properties of hyperbolic times.

\begin{proposition}[Proposition 2.6 in \cite{AA03}]
  \label{pr:backcontraction}
  There is $\de_1=\de_1(f_0)>0$ such that for every small
  enough $\ep>0$, if $n$ is $(\al,\delta)$-hyperbolic time
  for $(\omega,x)\in \supp(\theta^\ep)\times M$, then there
  is a neighborhood $V_n(\omega,x)$ of $x$ in $M$ such that
\begin{enumerate}
\item $f_{\omega}^n$ maps $V_n(\omega,x)$ diffeomorphically onto the ball
of radius $\de_1$ around $f_{\omega}^n(x)$;
\item for every $1\leq k\leq n$ and $y,z\in V_k(\omega,x)$
\[
\dist(f_{\omega}^{n-k}(y),f_{\omega}^{n-k}(z))
\leq
\al^{k/2}\dist(f_{\omega}^{n}(y),f_{\omega}^{n}(z)).
\]
\end{enumerate}
\end{proposition}

The uniform value of $\de_1$ in
Proposition~\ref{pr:backcontraction} is the crucial point to
get a uniform random generating partition. Indeed, let
$B_1,\dots, B_k$ be a finite open cover of $M$ by
$\de_1/2$-balls and let us take $\xi$ to be the partition
induced by this cover, i.e.
\[
\xi=\{ B_1, M\setminus B_1\} \vee \dots \vee \{ B_k, M\setminus B_k\}.
\]

\begin{lemma}[Lemma 6.6 in \cite{krerley2004}]
  \label{le:urandgenerating}
  If for a stationary measure $\mu$ we have that
  $\theta^\ep\times \mu$-almost all
  $(\omega,x)\in\Omega\times M$ have infinitely many
  $(\alpha,\de)$-hyperbolic times, then
\[
\lim_{k\to\infty}
  \diam( \vee_{j=0}^{k-1} (f_\omega^j)^{-1}\xi(x))= 0
\quad\mbox{for}\quad
\theta^\ep\times\mu-\mbox{almost every   } (\omega,x).
\]
\end{lemma}

By standard arguments, this ensures that if a stationary
measure $\mu$ is non-uniformly expanding and
$\theta^\ep\times\mu$-almost all random orbits have slow
approximation to $\cC$, then $\xi$ is a random generating
partition for $\mu$, see e.g.  Lemma 6.7 in
\cite{krerley2004}.

\subsubsection{Continuous variation of physical measures}
\label{sec:cont-vari-phys}

Here we prove Theorem \ref{thm.contNUE}.

We start by observing that since we are assuming that each
$f\in\cW$ has a physical measure $\mu_f$ satisfying
$m(M\setminus B(\mu_f))=0$, then by the proof of
Theorem~\ref{t.sts} we see that every weak$^*$ accumulation
point $\mu_0$ of the stationary measures $(\mu^{\ep_k})$
must equal $\mu_f$, since $\mu_0\ll m$. Thus $\mu_f$ is
stochastically stable for the random perturbations $(\hat
f,\hat \theta)$ we are considering.

In addition, since the exponent bound $c$ is uniform in
$\cW$, by the arguments in
Subsection~\ref{sec:unif-rand-gener} there is a uniform
generating partition $\xi$ for all $(f,\mu_f)$ with
$f\in\cW$.  Moreover for every $f\in\cW$ each $\mu_f$ is an
equilibrium state for $\phi_f=-\log|\det Df|$, $P_f(\phi_f)=0$ and
$\Phi:\cW\subset C^2(M,M)\to C^0(M,\bR),  f\mapsto \phi_f$ is
continuous. Then by Theorem~\ref{thm:contpressure} if we
take any sequence $f_k\in\cW$ converging to $f_0\in\cW$ when
$k\to\infty$, we know that every weak$^*$ accumulation
point $\mu_0$ of $(\mu_{f_k})_k$
satisfies~\eqref{eq:eqstate>0} with
$\phi_0=\phi_{f_0}$. Hence $\mu_0\ll m$ and by uniqueness of
the physical measure of $f_0$ we get $\mu_0=\mu_{f_0}$.

This finishes the proof of Theorem~\ref{thm.contNUE}.


\subsection{Local diffeomorphisms}
\label{sec:local-diff-1}

Here we prove the results in Subsection~\ref{sec:local-diff}.

In~\cite{AA03} it was shown that for the random
perturbations provided by Theorem~\ref{t.finite2} the
stochastic stability of $f_0$ implies non-uniform expansion
for random orbits. Here we prove the converse without the
extra technical condition used in~\cite{AA03}.

Let $(\hat f,\hat \theta)$ be a family of $C^2$ local
diffeomorphisms and of probability measures defining a
random perturbation of $f_0$ such that $f_0$ is non-uniformly
expanding and is non-uniformly expanding for random orbits.

First we observe that by Remark~\ref{rmk:nocritical} we may
use the results in Subsection~\ref{sec:unif-rand-gener}
also for local diffeomorphisms. Hence we can assume that
there exists a uniform random generating partition. We can
also use Theorem~\ref{t.sts} to conclude that any weak$^*$
accumulation point $\mu$ of stationary measure when $\ep\to0$  is
absolutely continuous.

Using now Theorem~\ref{thm:abv-Critical} we see that $\mu$
is a linear convex combination of stationary measures. This
means that $f_0$ is stochastically stable whenever it is
non-uniformly expanding for random orbits and ends the proof
of Theorem~\ref{thm.stc}.

\medskip

Finally, to prove Theorem~\ref{thm:semiconthtop} we just
have to use Theorem~\ref{thm:semicontpressure} with
$\phi_f\equiv0$ for all $f\in\cU$. This can be done since
uniform generating partitions exists for maximal entropy
measures.

\begin{proposition}[Lemma 4.8 in \cite{krerley2003}]
  \label{pr:hyptimesequilibrium}
There is $\de>0$ satisfying, for $f\in\cU$ (as defined in
Subsection~\ref{sec:equil-stat-potent}) and every equilibrium
state $\mu_\phi$ for a low-variation potential $\phi$ (as
given by Theorem~\ref{thm:krerley}),  that
$\mu_\phi$-almost every point $x\in M$ has infinitely many
$\de$-hyperbolic times.
\end{proposition}
This proposition together with
Proposition~\ref{pr:backcontraction} and
Lemma~\ref{le:urandgenerating} ensure the existence of a
fixed generating partition for every $f\in\cU$.

This concludes the proof of Theorem~\ref{thm:semiconthtop}.


  \bibliographystyle{abbrv}


\medskip

\end{document}